\documentclass[11pt,a4paper]{article}
\usepackage[cp1251]{inputenc}
\usepackage{amsfonts}
\usepackage{amssymb}
\usepackage{amsthm}
\usepackage{amsmath}
\usepackage{mathtext}
\usepackage{longtable}
\usepackage[russian]{babel}
\usepackage{color}

\title{Bounds for convergence rate in laws of large numbers\\ for mixed Poisson random sums\thanks{Supported by Russian
Science Foundation, project 18-11-00155.}}
\author{Victor Korolev\thanks{Faculty of Computational Mathematics
and Cybernetics, Moscow State University, Moscow, Russia; Institute
of Informatics Problems of the Federal Research Center ``Computer
Science and Control'',  Russian Academy of Sciences, Moscow, Russia.
E-mail: vkorolev@cs.msu.ru}, Alexander Zeifman\thanks{Vologda State
University, Vologda, Russia;  Institute of Informatics Problems of
the Federal Research Center ``Computer Science and Control'',
Russian Academy of Sciences, Moscow, Russia; Vologda Research Center
of the Russian Academy of Sciences, Vologda, Russia.  E-mail:
a$\_$zeifman@mail.ru}}
\date{}

\newcommand{\eqd}{\stackrel{d}{=}}

\renewcommand{\le}{\leqslant}
\renewcommand{\ge}{\geqslant}
\newcommand{\nb}{N\!\!\!B}

\begin{document}

\sloppy

\maketitle



\begin{center} {\bf Abstract} \end{center}

In the paper, upper bounds for the rate of convergence in laws of large numbers for mixed Poisson random sums are constructed. As a measure of the distance between the limit and pre-limit laws, the Zolotarev $\zeta$-metric is used. The obtained results extend the known convergence rate estimates for geometric random sums (in the famous R{\'e}nyi theorem) to a considerably wider class of random indices with mixed Poisson distributions including, e. g., those with the (generalized) negative binomial distribution.
\\[5pt]
\noindent {\bf Keywords:} law of large numbers; convergence rate; Zolotarev $\zeta$-metric; Poisson random sum; mixed Poisson distribution; geometric random sum; generalized negative binomial distribution\\[5pt]
\noindent {\bf AMS 2000 subject classification:} 60F05, 60G50,
60G55, 62E20, 62G30


\section{Introduction}


As far ago as in the 1950-s, being interested in modeling rare events, A. R{\'e}nyi studied rarefaction of renewal point processes and proved his famous theorem on convergence of rarefied renewal processes to the Poisson process \cite{Renyi1956, Renyi1964}. He considered a simple rarefaction, that is, each point either remains as it is with probability $p$ or is removed with probability $1-p$ independently of other points. A. R{\'e}nyi showed that, after infinitely many iterations of rarefaction accompanied by an appropriate contraction of time, the resulting process is Poisson characterized by the exponential distribution of time intervals between successive points. The key point was that the interval between successive points in the rarefied process is a sum of a random number of independent random variables in which the number of summands has the geometric distribution and is independent of the summands, that is, a geometric sum. Geometric sums appeared to be important mathematical models in risk theory and insurance \cite{Kalashnikov1997, Grandell2000, Grandell2002}, reliability theory \cite{Kalashnikov1997, Bon1999, Bon2002} and other fields. An outburst of interest to {\it analytic} and asymptotic properties of geometric sums was initiated by the publication of the path-breaking paper of L. B. Klebanov and his colleagues \cite{Klebanovetal1984} introducing the notion of geometric infinite divisibility.

The R{\'e}nyi theorem states that the distribution of a geometric sum normalized by its expectation converges to the exponential law as the expectation of the sum infinitely increases. The normalization of a sum by its expectation is typical for laws of large numbers. Therefore, the R{\'e}nyi theorem can be regarded as the law of large numbers for geometric sums. A general law of large numbers for random sums of independent identically distributed (i.i.d.) random variables (r.v.'s) was proved in \cite{Korolev1994}. It was demonstrated there that the distribution of a random sum normalized by its expectation converges to some distribution, if and only if the distribution of the random index (the number of summands) converges to the same distribution (up to a scale parameter) under the same normalization.

First studies of the rate of convergence in the law of large numbers for random sums were quite naturally concentrated around the R{\'e}nyi theorem. The first result was due to A. D. Solovyev \cite{Solovyev1971}. V. V. Kalashnikov and S. Y. Vsekhsvyatskii \cite{KV1985} were the first to obtain bounds of convergence rate in the R{\'e}nyi theorem in terms of the Zolotarev ideal $\zeta$-metric. The best known bounds are due to M. Brown \cite{Brown1990} (for non-negative summands) and I. Shevtsova and M. Tselishchev \cite{ShevtsovaTselishchev2020} (for the general case).

A natural direction of development of these studies is to extend the results to more general settings concerning the distribution of the number of summands. A direct way is to consider negative binomial random sums with geometric sums being a particular case. Negative binomial random sums play an important role in modeling real phenomena, see, e. g., \cite{SheejaKumar2017} for financial applications and \cite{KorolevGorsheninBelyaev2019} for climatological applications. Limit theorems for negative binomial sums were considered in many papers, see, e. g., \cite{BeningKorolev2005, KorolevDAN,  SchluterTrede2016} and the references therein. Many papers dealt with bounds of convergence rate in analogs of central limit theorem for negative binomial sums, see, e. g., \cite{KorolevShevtsovaSAJ, KorolevDorofeevaLMJ}, where the uniform (Kolmogorov) distance was considered or \cite{Sunklodas2015} where bounds were constructed in terms of the $\zeta$-metric. As regards approximation of the distribution of negative binomial sums by the gamma distribution, we should mention the paper \cite{Gavrilenko2017} where bounds were constructed for the accuracy of gamma-approximation to the negative binomial distribution with small shape parameter $r$, and the paper \cite{TranLocHung2018} presenting the bounds for the rate of convergence of the distribution of negative binomial sums to the gamma distribution in terms of the Trotter distance. Strange as it may seem, in this setting there are still no bounds in terms of the convenient $\zeta$-metric (at least as far as the authors know).

The present paper fills this gap. Moreover, here we present a technique that makes it possible to construct the bounds in terms of the $\zeta$-metric for the rate of convergence in the laws of large numbers for random sums with mixed Poisson-distributed number of summands covering the geometric, negative binomial and many other important and useful cases including the so-called generalized negative binomial sums (see, e. g., \cite{KorolevZeifman2019, KorolevGorshenin2020}). It should be noted that the results presented below do not assume any restrictions on the range of the summands. It is only required that the second moment of summands exists and the expectation differs from zero.

\section{Mathematical preliminaries}

In the paper, conventional notation is used. The symbols $\eqd$ and $\Longrightarrow$ denote the coincidence of distributions and convergence in distribution, respectively.

In what follows, for brevity and convenience, the results will be presented in terms of r.v.'s with the corresponding distributions. It will be assumed that all the r.v.'s are defined on the same probability space $(\Omega,\,\mathfrak{F},\,{\sf P})$.

Let $s>0$. There exist a unique representation of the number $s$ as $s=m+\alpha$ where $m$ is an integer and $0<\alpha\le1$. By $\mathcal{F}_s$ we denote the set of all real-valued bounded functions $f$ on $\mathbb{R}$ that are $m$ times differentiable and
\begin{equation}\label{1}
|f^{(m)}(x)-f^{(m)}(y)|\le |x-y|^{\alpha}.
\end{equation}
Let $X$ and $Y$ be two r.v.'s whose distribution functions will be denoted $F_X(x)$ and $F_Y(x)$, respectively. In \cite{Zolotarev1976} V. M. Zolotarev introduced a convenient {\it ideal} $\zeta$-metric $\zeta_s(X,\,Y)\equiv\zeta_s(F_X,\,F_Y)$ in the space of probability distributions by the equality
\begin{equation}\label{2}
\zeta_s(X,\,Y)=\sup\big\{\big|{\sf E}\big(f(X)-f(Y)\big)\big|:\,f\in\mathcal{F}_s\big\},
\end{equation}
also see \cite{Zolotarev1977} and \cite{Zolotarev1997}, p. 44. In particular,
$$
\zeta_1(X,\,Y)=\int_{\mathbb{R}}|F_X(x)-F_Y(x)|dx.
$$

For convenience, the properties of $\zeta$-metrics will be formulated as lemmas.

\smallskip

{\sc Lemma 1}. {\it Let $c>0$. Then}
$$
\zeta_s(cX,\,cY)=c^s\zeta_s(X,\,Y).
$$

\smallskip

{\sc Lemma 2}. {\it Let $X,Y,Z$ be r.v.'s such that $Z$ is independent of both $X$ and $Y$. Then}
$$
\zeta_s(X+Z,\,Y+Z)=\zeta_s(X,\,Y).
$$

\smallskip

The proofs of these statements can be found in \cite{Zolotarev1997}, see Theorem 1.4.2 there. Lemma 1 states the property of {\it homogeneity} of order $s$ of the metric $\zeta_s$ and Lemma 2 states its property of {\it regularity}. The proof of regularity of the metric $\zeta_s$ given in \cite{Zolotarev1997} can be easily extended from convolutions, i. e. location or shift mixtures, to arbitrary mixtures.

\smallskip

{\sc Lemma 3}. {\it Let $X$ and $Y$ be r.v.'s with the distribution functions $F_X(x)$ and $F_Y(x)$, respectively. Let $F_X(x;z)$ and $F_Y(x;z)$ be functions of two variables $x,z\in\mathbb{R}$ and such that they are measurable with respect to $z$ for each fixed $x$ and are distribution functions as functions of $x$ for each fixed $z$. Let $H(z)$ be the distribution function of a r.v. $Z$ such that for $f\in\mathcal{F}_s$
$$
{\sf E}\big[f(X)|Z=z\big]=\int_{\mathbb{R}}f(x)d_xF_X(x;z),\ \ \ \  {\sf E}\big[f(Y)|Z=z\big]=\int_{\mathbb{R}}f(x)d_xF_Y(x;z).
$$
Then}
$$
\zeta_s(F_X,\,F_Y)\le\int_{\mathbb{R}}\zeta_s\big(F_X(\,\cdot\,;z),\,F_Y(\,\cdot\,;z)\big)dH(z).
$$

\smallskip

{\sc Proof}. For any $f\in\mathcal{F}_s$ by the formula of total expectation we have
$$
{\sf E}f(X)={\sf E}{\sf E}\big[f(X)|Z\big]=\int_{\mathbb{R}}{\sf E}\big[f(X)|Z=z\big]dH(z)=\int_{\mathbb{R}}\int_{\mathbb{R}}f(x)d_xF_X(x;z)dH(z),
$$
$$
{\sf E}f(Y)={\sf E}{\sf E}\big[f(Y)|Z\big]=\int_{\mathbb{R}}{\sf E}\big[f(Y)|Z=z\big]dH(z)=\int_{\mathbb{R}}\int_{\mathbb{R}}f(x)d_xF_Y(x;z)dH(z),
$$
so that {\it for each} $f\in\mathcal{F}_s$
$$
\big|{\sf E}f(X)-{\sf E}f(Y)\big|=\bigg|\int_{\mathbb{R}}\bigg[\int_{\mathbb{R}}f(x)d_xF_X(x;z)-\int_{\mathbb{R}}f(x)d_xF_Y(x;z)\bigg]dH(z)\bigg|\le
$$
$$
\le\int_{\mathbb{R}}\bigg|\int_{\mathbb{R}}f(x)d_xF_X(x;z)-\int_{\mathbb{R}}f(x)d_xF_Y(x;z)\bigg|dH(z)\le
$$
$$
\le\int_{\mathbb{R}}\sup\bigg\{\bigg|\int_{\mathbb{R}}f(x)d_xF_X(x;z)-\int_{\mathbb{R}}f(x)d_xF_Y(x;z)\bigg|:\,f\in\mathcal{F}_s\bigg\}dH(z)=
$$
$$
=\int_{\mathbb{R}}\zeta_s\big(F_X(\,\cdot\,,z),\,F_Y(\,\cdot\,,z)\big)dH(z).
$$
The remark that the inequality
$$
\big|{\sf E}f(X)-{\sf E}f(Y)\big|\le \int_{\mathbb{R}}\zeta_s\big(F_X(\,\cdot\,,z),\,F_Y(\,\cdot\,,z)\big)dH(z)
$$
just established holds uniformly in $f\in\mathcal{F}_s$ completes the proof.

\smallskip

{\sc Lemma 4}. {\it Let $s>0$, $s=m+\alpha$ with $0<\alpha\le1$. Let $X$ and $Y$ be r.v.'s such that ${\sf E}|X|^s$ and
${\sf E}|Y|^s$ are finite and, moreover, ${\sf E}X^k={\sf E}Y^k$ for $k=1,\ldots,m$. Then}
$$
\zeta_s(X,\,Y)\le \frac{\Gamma(1+\alpha)}{\Gamma(1+s)}\cdot\big[{\sf E}|X|^s+{\sf E}|Y|^s\big].
$$

\smallskip

{\sc Proof}. This statement is a simple consequence of Theorem 1.7 (c) in \cite{MattnerShevtsova2019}. A similar inequality (however, without the term $\Gamma(1+\alpha)\le1$) can be also deduced from relation (1.5.41) in \cite{Zolotarev1997} (see Theorem 1.5.7 on p. 78 there).

\smallskip

Consider i.i.d. r.v.'s $X_1,X_2,\ldots$. Let $\lambda>0$ and $N(\lambda)$ be the r.v. with the Poisson distribution
$$
{\sf P}(N(\lambda)=k)=e^{-\lambda}\frac{\lambda^k}{k!},\ \ \ \ \ k=0,1,2,\ldots
$$
Assume that for each $\lambda>0$ the r.v. $N(\lambda)$ is independent of $X_1,X_2,\ldots$. It is well known that if $a\equiv{\sf E}X_1$ is finite, then
\begin{equation}\label{4}
\frac{1}{\lambda}\sum\nolimits_{j=1}^{N(\lambda)}X_j\Longrightarrow a
\end{equation}
as $\lambda\to\infty$ (for definiteness, the sum $\sum\nolimits_{j=1}^0$ is assumed to be equal to zero). The following lemma establishes the rate of convergence in \eqref{4}.

\smallskip

{\sc Lemma 5}. {\it In addition to the above assumptions, let $a\neq 0$ and $\sigma^2\equiv {\sf D}X_1$ be finite. Then for $1\le s\le2$ we have}
\begin{equation}\label{5}
\zeta_s\bigg(\frac{1}{a\lambda}\sum\nolimits_{j=1}^{N(\lambda)}X_j,\,1\bigg)\le\frac{\Gamma(1+\alpha)}{\Gamma(1+s)}\cdot\Big(\frac{a^2+\sigma^2}{\lambda a^2}\Big)^{s/2}.
\end{equation}

\smallskip

{\sc Proof}. First, note that
\begin{equation}\label{110}
{\sf D}\sum\nolimits_{j=1}^{N(\lambda)}X_j=\lambda(a^2+\sigma^2).
\end{equation}
Successively applying Lemma 2, Lemma 4, Lemma 1, the Jensen inequality for $1\le s\le2$ and \eqref{110}, we obtain
$$
\zeta_s\bigg(\frac{1}{a\lambda}\sum\nolimits_{j=1}^{N(\lambda)}X_j,\,1\bigg)=
\zeta_s\bigg(\frac{1}{a\lambda}\sum\nolimits_{j=1}^{N(\lambda)}X_j-1,\,0\bigg)\le\frac{\Gamma(1+\alpha)}{(\lambda a)^s\Gamma(1+s)}
{\sf E}\bigg|\sum\nolimits_{j=1}^{N(\lambda)}X_j-\lambda a\bigg|^s\le
$$
$$
\le\frac{\Gamma(1+\alpha)}{(\lambda a)^s\Gamma(1+s)}\bigg({\sf E}\bigg[\sum\nolimits_{j=1}^{N(\lambda)}X_j-\lambda a\bigg]^2\bigg)^{s/2}\!=
\frac{\Gamma(1+\alpha)}{(\lambda a)^s\Gamma(1+s)}\bigg({\sf D}\sum\nolimits_{j=1}^{N(\lambda)}X_j\bigg)^{s/2}=
$$
\begin{equation}\label{6}
=\frac{\Gamma(1+\alpha)}{\Gamma(1+s)}\cdot\Big(\frac{a^2+\sigma^2}{\lambda a^2}\Big)^{s/2}.
\end{equation}
The lemma is proved.

\smallskip

If $s=2$, then \eqref{5} turns into
$$
\zeta_2\bigg(\frac{1}{a\lambda}\sum\nolimits_{j=1}^{N(\lambda)}X_j,\,1\bigg)\le\frac{1}{2\lambda}\Big(1+\frac{\sigma^2}{a^2}\Big).
$$

\section{Main results}

Consider a sequence of positive r.v.'s $\Lambda_1,\Lambda_2,\ldots$. Let $N(t)$ be the standard Poisson process (the Poisson process with unit intensity) independent of the sequence $\Lambda_1,\Lambda_2,\ldots$. Let for each $n\in\mathbb{N}$ the r.v. $N_n$ be defined as
$$
N_n=N(\Lambda_n).
$$
The r.v. $N_n$ so defined has a {\it mixed Poisson distribution}:
$$
{\sf P}(N_n=k)=\frac{1}{k!}\int_{0}^{\infty}e^{-\lambda}\lambda^kd{\sf P}(\Lambda_n<\lambda),\ \ \ k=0,1,2,\ldots.
$$
The class of mixed Poisson distributions is very wide and contains, for example, the geometric distribution \cite{Kalashnikov1997} (in that case the distribution of $\Lambda_n$ is exponential), the negative binomial distribution \cite{GreenwoodYule1920} (with the gamma distribution of $\Lambda_n$), the generalized negative binomial distributions \cite{KorolevZeifman2019} (with the generalized gamma-distributed $\Lambda_n$), the Sichel distribution \cite{Holla1967, Sichel1971} (with $\Lambda_n$ having the generalized inverse Gaussian distribution) and many other types, see Example 3 below.

Again consider a sequence $X_1,X_2,\ldots$ of i.i.d. r.v.'s with finite ${\sf E}X_1\equiv a\neq 0$. Assume that for each $n\in\mathbb{N}$ the r.v. $N_n$ is independent of the sequence $X_1,X_2,\ldots$ We will concentrate our attention on the asymptotic behavior of {\it mixed Poisson random sums}
\begin{equation}\label{7}
S_n=\sum\nolimits_{j=1}^{N_n}X_j.
\end{equation}
Let $\{m_n\}_{n\ge1}$ be an infinitely increasing sequence of positive numbers. It is known that a r.v. $\Lambda$ such that
\begin{equation}\label{8}
\frac{S_n}{m_n}\Longrightarrow a\Lambda \ \ \ \ (n\to\infty)
\end{equation}
exists if and only if
\begin{equation}\label{9}
\frac{\Lambda_n}{m_n}\Longrightarrow \Lambda \ \ \ \ (n\to\infty)
\end{equation}
(see, e. g., \cite{Korolev1994, Korolev1998}). Relation \eqref{8} can be referred to as the law of large numbers for (mixed Poisson) random sums. Unlike the classical laws of large numbers for sums of a non-random number of r.v.'s, for random sums the limit can be random and coincides with the limit for mixing r.v.'s under the same normalization.

Consider the rate of convergence in \eqref{8}. For convenience and better vividness, we will sometimes use the full notation $N(\Lambda_n)$ for $N_n$.

\smallskip

{\sc Theorem 1}. {\it In addition to the above assumptions, let $\sigma^2\equiv {\sf D}X_1$ be finite. Then for $1\le s\le2$ we have}
\begin{equation}\label{10}
\zeta_s\bigg(\frac{1}{am_n}\sum\nolimits_{j=1}^{N(\Lambda_n)}X_j,\,\Lambda\bigg)\le \frac{{\sf E}\Lambda_n^{s/2}}{m_n^s}\cdot\frac{\Gamma(1+\alpha)}{\Gamma(1+s)}\cdot\Big(1+\frac{\sigma^2}{a^2}\Big)^{s/2}+
\zeta_s\Big(\frac{\Lambda_n}{m_n},\,\Lambda\Big).
\end{equation}

\smallskip

{\sc Proof}. By the triangle inequality we have
\begin{equation}\label{11}
\zeta_s\bigg(\frac{1}{am_n}\sum\nolimits_{j=1}^{N(\Lambda_n)}X_j,\,\Lambda\bigg)\le
\zeta_s\bigg(\frac{1}{am_n}\sum\nolimits_{j=1}^{N(\Lambda_n)}X_j,\,\frac{\Lambda_n}{m_n}\bigg)+\zeta_s\Big(\frac{\Lambda_n}{m_n},\,\Lambda\Big).
\end{equation}
In order to estimate the first term on the right-hand side of \eqref{11} we successively apply Lemma 3, Lemma 2, Lemma 4, Lemma 1, the Jensen inequality for $1\le s\le2$ and \eqref{110} and obtain the following chain of relations:
$$
\zeta_s\bigg(\frac{1}{am_n}\sum\nolimits_{j=1}^{N(\Lambda_n)}X_j,\,\frac{\Lambda_n}{m_n}\bigg)\le \int_{0}^{\infty}\zeta_s\bigg(\frac{1}{am_n}\sum\nolimits_{j=1}^{N(\lambda)}X_j,\,\frac{\lambda}{m_n}\bigg)d{\sf P}(\Lambda_n<\lambda)=
$$
$$
=\int_{0}^{\infty}\zeta_s\bigg(\frac{1}{am_n}\sum\nolimits_{j=1}^{N(\lambda)}X_j-\frac{\lambda}{m_n},\,0\bigg)d{\sf P}(\Lambda_n<\lambda)=
\frac{1}{(am_n)^s}\int_{0}^{\infty}\zeta_s\bigg(\sum\nolimits_{j=1}^{N(\lambda)}X_j-\lambda a,\,0\bigg)d{\sf P}(\Lambda_n<\lambda)\le
$$
$$
\le\frac{\Gamma(1+\alpha)}{\Gamma(1+s)(am_n)^s}\int_{0}^{\infty}{\sf E}\bigg|\sum\nolimits_{j=1}^{N(\lambda)}X_j-\lambda a\bigg|^sd{\sf P}(\Lambda_n<\lambda)\le
$$
$$
\le\frac{\Gamma(1+\alpha)}{\Gamma(1+s)(am_n)^s}\int_{0}^{\infty}\bigg({\sf D}\sum\nolimits_{j=1}^{N(\lambda)}X_j\bigg)^{s/2}d{\sf P}(\Lambda_n<\lambda)=
$$
$$
= \frac{\Gamma(1+\alpha)}{\Gamma(1+s)(am_n)^s}\Big(1+\frac{\sigma^2}{a^2}\Big)^{s/2}\int_{0}^{\infty}\lambda^{s/2}d{\sf P}(\Lambda_n<\lambda)=
\frac{{\sf E}\Lambda_n^{s/2}}{m_n^s}\frac{\Gamma(1+\alpha)}{\Gamma(1+s)}\Big(1+\frac{\sigma^2}{a^2}\Big)^{s/2}.
$$
The theorem is proved.

\smallskip

{\sc Remark 1}. If $m_n={\sf E}\Lambda_n$, then by virtue of the Jensen inequality instead of \eqref{10} we have
\begin{equation}\label{111}
\zeta_s\bigg(\frac{1}{am_n}\sum\nolimits_{j=1}^{N(\Lambda_n)}X_j,\,\Lambda\bigg)\le
\frac{\Gamma(1+\alpha)}{m_n^{s/2}\Gamma(1+s)}\Big(1+\frac{\sigma^2}{a^2}\Big)^{s/2}+\zeta_s\Big(\frac{\Lambda_n}{m_n},\,\Lambda\Big)
\end{equation}
since $s/2\le1$ for the specified range of $s$. In particular,
\begin{equation}\label{112}
\zeta_2\bigg(\frac{1}{am_n}\sum\nolimits_{j=1}^{N(\Lambda_n)}X_j,\,\Lambda\bigg)\le
\frac{1}{2m_n}\Big(1+\frac{\sigma^2}{a^2}\Big)+\zeta_2\Big(\frac{\Lambda_n}{m_n},\,\Lambda\Big)
\end{equation}

\smallskip

{\sc Remark 2}. Since ${\sf E}X_1^2=a^2+\sigma^2$, we have
$$
1+\frac{\sigma^2}{a^2}=\frac{{\sf E}X_1^2}{({\sf E}X_1)^2}.
$$

\smallskip

{\sc Corollary 1}. {\it Under the conditions of Theorem 1 we have}
$$
\zeta_2\bigg(\frac{1}{am_n}\sum\nolimits_{j=1}^{N(\Lambda_n)}X_j,\,\Lambda\bigg)\le
\frac{{\sf E}\Lambda_n}{2m_n^2}\Big(1+\frac{\sigma^2}{a^2}\Big)+\zeta_2\Big(\frac{\Lambda_n}{m_n},\,\Lambda\Big).
$$

\smallskip

Now consider the case where, as $n$ varies, the r.v.'s $\Lambda_n$ preserve the same distribution type with the only change of scale factor.

\smallskip

{\sc Corollary 2}. {\it Assume that, in addition to the conditions of Theorem 1, $\Lambda_n=m_n\Lambda$. Then}
\begin{equation}\label{13}
\zeta_s\bigg(\frac{1}{am_n}\sum\nolimits_{j=1}^{N(m_n\Lambda)}\!\!X_j,\,\Lambda\bigg)\le
\frac{{\sf E}\Lambda^{s/2}}{m_n^{s/2}}\frac{\Gamma(1+\alpha)}{\Gamma(1+s)}\Big(1+\frac{\sigma^2}{a^2}\Big)^{s/2}.
\end{equation}
{\it In particular,}
$$
\zeta_2\bigg(\frac{1}{am_n}\sum\nolimits_{j=1}^{N(m_n\Lambda)}\!\!X_j,\,\Lambda\bigg)\le
\frac{{\sf E}\Lambda}{2m_n}\Big(1+\frac{\sigma^2}{a^2}\Big).
$$

\smallskip

{\sc Example 1}. Let $\Lambda$ have the standard exponential distribution and $\Lambda_n=m_n\Lambda$. Then $N_n=N(m_n\Lambda)$ has the geometric distribution with parameter $p_n=(1+m_n)^{-1}$. In this case $m_n=(1-p_n)/p_n$, ${\sf E}\Lambda=1$ and from \eqref{111} we obtain the following convergence rate bound in the R{\'e}nyi theorem:
$$
\zeta_s\bigg(\frac{p_n}{1-p_n}\sum\nolimits_{j=1}^{N_n}X_j,\,\Lambda\bigg)\le \frac{\Gamma(1+\alpha)}{\Gamma(1+s)}\bigg[\frac{p_n}{1-p_n}\cdot\frac{{\sf E}X_1^2}{({\sf E}X_1)^2}\bigg]^{s/2}, \ \ \ 1\le s\le2.
$$
In particular,
$$
\zeta_2\bigg(\frac{p_n}{1-p_n}\sum\nolimits_{j=1}^{N_n}X_j,\,\Lambda\bigg)\le \frac{p_n}{2(1-p_n)}\cdot\frac{{\sf E}X_1^2}{({\sf E}X_1)^2}
$$
(sf. \cite{Kalashnikov1997} and \cite{ShevtsovaTselishchev2020}).

\smallskip

{\sc Example 2}. Let $\Lambda_n$ have the gamma distribution with shape parameter $r>0$ and scale parameter $\mu_n=\mu/n$, $\mu>0$. Then $\Lambda_n=G_{r,\mu/n}\eqd nG_{r,\mu}$ where $G_{r,\mu}$ has the gamma distribution with parameters $r$ and $\mu$ corresponding to the probability density function
$$
g(x;\,r,\mu)=\frac{\mu^r}{\Gamma(r)}x^{r-1}e^{-\mu x},\ \ \ x\ge0.
$$
Then $N_n=N(nG_{r,\mu})$ has the negative binomial distribution with parameters $r$ and $p_n=\mu/(n+\mu)$:
$$
{\sf P}(N_n=k)=\frac{\mu^rn^k}{k!\Gamma(r)}\int_{0}^{\infty}e^{-\lambda(n+\mu)}\lambda^{k+r-1}d\lambda=
\frac{\mu^rn^k}{k!\Gamma(r)(n+\mu)^{k+r}}\int_{0}^{\infty}e^{-\lambda}\lambda^{k+r-1}d\lambda=
$$
$$
=\frac{\Gamma(k+r)}{k!\Gamma(r)}\Big(\frac{\mu}{n+\mu}\Big)^r\Big(1-\frac{\mu}{n+\mu}\Big)^k,\ \ \ \ k=0,1,2,\ldots
$$
We have $m_n={\sf E}G_{r,\mu/n}=n{\sf E}(G_{r,\mu})=nr/\mu$ so that for each $n\in\mathbb{N}$
$$
\frac{\Lambda_n}{m_n}=\frac{nG_{r,\mu}}{m_n}=\frac{\mu}{r}G_{r,\mu}\eqd G_{r,r},
$$
that is, the limit distribution is gamma with shape and scale parameters equal to $r$: $\Lambda\eqd G_{r,r}$. Therefore, for $1\le s\le2$ due to Corollary 2 the following bound holds:
\begin{equation}\label{12}
\zeta_s\bigg(\frac{\mu}{anr}\sum\nolimits_{j=1}^{N(nG_{r,\mu})}X_j,\,G_{r,r}\bigg)\le\frac{\Gamma(1+\alpha)}{\Gamma(1+s)}\bigg(\frac{\mu}{nr}\cdot\frac{{\sf E}X_1^2}{({\sf E}X_1)^2}\bigg)^{s/2}.
\end{equation}
In particular,
\begin{equation}\label{13}
\zeta_2\bigg(\frac{\mu}{anr}\sum\nolimits_{j=1}^{N(nG_{r,\mu})}X_j,\,G_{r,r}\bigg)\le\frac{\mu}{2nr}\cdot\frac{{\sf E}X_1^2}{({\sf E}X_1)^2}.
\end{equation}
To emphasize the fact that $N(nG_{r,\mu})$ has the negative binomial distribution with parameters $r$ and $p_n=\mu/(n+\mu)$, introduce an alternative notation for this r.v.: $N(nG_{r,\mu})=\nb_{r,p_n}$. In these terms \eqref{12} and \eqref{13} can be rewritten as
\begin{equation}\label{14}
\zeta_s\bigg(\frac{p_n}{ar(1-p_n)}\sum\nolimits_{j=1}^{\nb_{r,p_n}}X_j,\,G_{r,r}\bigg)\le
\frac{\Gamma(1+\alpha)}{\Gamma(1+s)}\bigg[\frac{p_n}{(1-p_n)r}\cdot\frac{{\sf E}X_1^2}{({\sf E}X_1)^2}\bigg]^{s/2}
\end{equation}
and
\begin{equation}\label{15}
\zeta_2\bigg(\frac{p_n}{ar(1-p_n)}\sum\nolimits_{j=1}^{\nb_{r,p_n}}X_j,\,G_{r,r}\bigg)\le
\frac{p_n}{2r(1-p_n)}\cdot\frac{{\sf E}X_1^2}{({\sf E}X_1)^2}.
\end{equation}

\smallskip

{\sc Example 3}. Let $\Lambda_n=G^*_{r,\alpha,\mu/n^{\alpha}}\eqd nG^*_{r,\alpha,\mu}=$ where the r.v. $G^*_{r,\alpha,\mu}$ has the generalized gamma distribution (GG distribution) corresponding to the density
$$
g^*(x;r,\alpha,\mu)=\frac{|\alpha|\mu^r}{\Gamma(r)}x^{\alpha
r-1}e^{-\mu x^{\alpha}},\ \ \ \ x\ge0,
$$
with $\alpha\in\mathbb{R}$, $\mu>0$, $r>0$.

The class of GG distributions was first described as a unitary family in 1962 by E. Stacy \cite{Stacy1962} as the class of probability distributions
simultaneously containing both Weibull and gamma distributions. The family of GG distributions contains practically all the most popular absolutely continuous distributions concentrated on the non-negative half-line. In particular, the family of GG distributions contains:

\begin{itemize}

\vspace{-3mm}\item[$\bullet$] the gamma distribution $(\alpha=1)$
and its special cases:

\begin{itemize}

\vspace{-3mm}\item[$\circ$] the exponential distribution
($\alpha=1$, $r =1$),

\vspace{-1mm}\item[$\circ$] the Erlang distribution ($\alpha=1$, $ r
\in\mathbb{N}$),

\vspace{-1mm}\item[$\circ$] the chi-square distribution ($\alpha=1$,
$\mu=\frac12$);

\end{itemize}

\vspace{-4mm}\item[$\bullet$] the Nakagami distribution
($\alpha=2$);

\vspace{-2mm}\item[$\bullet$] the half-normal (folded normal)
distribution (the distribution of the maximum of a standard Wiener
process on the interval $[0,1]$) ($\alpha=2$, $ r =\frac12$);

\vspace{-2mm}\item[$\bullet$] the Rayleigh distribution ($\alpha=2$,
$ r =1$);

\vspace{-2mm}\item[$\bullet$] the chi-distribution ($\alpha=2$,
$\mu=1/\sqrt{2}$);

\vspace{-2mm}\item[$\bullet$] the Maxwell distribution (the
distribution of the absolute values of the velocities of moleculas
in a dilute gas) ($\alpha=2$, $ r =\frac32$);

\vspace{-2mm}\item[$\bullet$] the Weibull--Gnedenko distribution
(the extreme value distribution of type III) ($r =1$, $\alpha>0$);

\vspace{-2mm}\item[$\bullet$] the (folded) exponential power
distribution ($\alpha>0$, $r=\frac{1}{\alpha}$);

\vspace{-2mm}\item[$\bullet$] the inverse gamma distribution
($\alpha=-1$) and its special case:

\begin{itemize}

\vspace{-3mm}\item[$\circ$] the L{\'e}vy distribution (the one-sided
stable distribution with the characteristic exponent $\frac12$ --
the distribution of the first hit time of the unit level by the
Brownian motion) ($\mu=-1$, $r =\frac12$);

\end{itemize}

\vspace{-2mm}\item[$\bullet$] the Fr{\'e}chet distribution (the
extreme value distribution of type II) ($r=1$, $\alpha<0$)

\end{itemize}

\vspace{-2mm} \noindent and other laws. The limit point of the class
of GG distributions is

\begin{itemize}

\vspace{-2mm}\item[$\bullet$] the log-normal distribution ($r
\to\infty$).

\end{itemize}

\vspace{-2mm}

There are dozens of papers dealing with the application of GG distributions as models of regularities observed in practice. Apparently, the popularity of GG distributions is due to that most of them can serve as adequate asymptotic approximations, since all the representatives of the class of GG distributions listed above appear as limit laws in various limit theorems of probability theory in rather simple limit schemes. In \cite{KorolevZeifman2019} a general limit theorem (an analog of the law of large numbers) for random sums of independent r.v.'s was formulated in which the GG distributions are limit laws.

It is easy to make sure that the moment of order $\delta>0$ of $G^*_{r,\alpha,\mu}$ has the form
$$
{\sf E}(G^*_{r,\alpha,\mu})^{\delta}=\frac{\Gamma(r+\frac{\delta}{\alpha})}{\mu^{\delta/\alpha}\Gamma(r)}.
$$
If we take $m_n=n$, then the limit r.v. $\Lambda=G^*_{r,\alpha,\mu}$ and from Corollary 2 it follows that for $1\le s\le2$
$$
\zeta_s\bigg(\frac{1}{na}\sum\nolimits_{j=1}^{N(nG^*_{r,\alpha,\mu})}\!X_j,\,G^*_{r,\alpha,\mu}\bigg)\le
\frac{({\sf E}X_1^2)^{s/2}}{n^{s/2}|{\sf E}X_1|^s}\cdot\frac{}{}\frac{\Gamma(1+\alpha)\Gamma(r+\frac{s}{2\alpha})}{\mu^{s/(2\alpha)}\Gamma(1+s)\Gamma(r)}.
$$
In particular,
$$
\zeta_2\bigg(\frac{1}{na}\sum\nolimits_{j=1}^{N(nG^*_{r,\alpha,\mu})}\!X_j,\,G^*_{r,\alpha,\mu}\bigg)\le
\frac{{\sf E}X_1^2}{2n({\sf E}X_1)^2}\cdot\frac{\Gamma(r+\frac{1}{\alpha})}{\mu^{1/\alpha}\Gamma(r)}.
$$

\renewcommand{\refname}{Acknowledgement}

\vspace{-0.3cm}

The authors have the pleasure to thank Irina Shevtsova for valuable bibliographical advices.

\smallskip

\renewcommand{\refname}{References}

\end{document}